\documentclass[11pt]{amsart}
\usepackage{amsxtra}
\usepackage{amssymb}
\theoremstyle{plain}

\newcommand{\cleqn}{\setcounter{equation}{0}}
\newcommand{\clth}{\setcounter{theorem}{0}}
\newcommand {\sectionnew}[1]{\section{#1}\cleqn\clth}

\newtheorem{theorem}{Theorem}[section]
\newtheorem{lemma}[theorem]{Lemma}
\newtheorem{definition-theorem}[theorem]{Definition-Theorem}
\newtheorem{proposition}[theorem]{Proposition}
\newtheorem{corollary}[theorem]{Corollary}
\newtheorem{definition}[theorem]{Definition}
\newtheorem{example}[theorem]{Example}
\newtheorem{remark}[theorem]{Remark}
\newtheorem{conjecture}[theorem]{Conjecture}
\newtheorem{notation}[theorem]{Notation}
\newcommand \bth[1] { \begin{theorem}\label{t#1} }
\newcommand \ble[1] { \begin{lemma}\label{l#1} }

\newcommand \bpr[1] { \begin{proposition}\label{p#1} }
\newcommand \bco[1] { \begin{corollary}\label{c#1} }
\newcommand \bde[1] { \begin{definition}\label{d#1}\rm }
\newcommand \bex[1] { \begin{example}\label{e#1}\rm }
\newcommand \bre[1] { \begin{remark}\label{r#1}\rm }
\newcommand \bcj[1] { \begin{conjecture}\label{j#1}\rm }

\newcommand \bnota[1] { \begin{notation}\label{n#1}\rm }
\renewcommand {\eth} { \end{theorem} }
\newcommand {\ele} { \end{lemma} }

\newcommand {\epr} { \end{proposition} }
\newcommand {\eco} { \end{corollary} }
\newcommand {\ede} { \end{definition} }
\newcommand {\eex} { \end{example} }
\newcommand {\ere} { \end{remark} }
\newcommand {\ecj} { \end{conjecture} }

\newcommand {\enota} { \end{notation} }
\newcommand \thref[1]{Theorem \ref{t#1}}
\newcommand \leref[1]{Lemma \ref{l#1}}
\newcommand \prref[1]{Proposition \ref{p#1}}

\newcommand \cjref[1]{Conjecture \ref{j#1}}

\newcommand \reref[1]{Remark \ref{r#1}}
\newcommand \lb[1]{\label{#1}}

\def \Cset {{\mathbb C}}
\def \KK {{\mathbb K}}
\def \Zset {{\mathbb Z}}

\def \Qset {{\mathbb Q}}

\def \B  {{\mathcal{B}}}
\def \II {{\mathcal{I}}}
\def \QQ {{\mathcal{Q}}}
\def \JJ {{\mathcal{J}}}

\def \UU {{\mathcal{U}}}
\def \RR {{\mathcal{R}}}

\def \LL {{\mathcal{L}}}
\def \TT {{\mathcal{T}}} 

\def \al {\alpha}

\def \la {\lambda}

\def \om {\omega}
\def \ga {\gamma}

\def \mt  {\mapsto}

\def \lha {\leftharpoonup}
\def \rha {\rightharpoonup}

\def \rcor {\rangle}
\def \lcor {\langle}

\def \ol {\overline}
\def \wt {\widetilde}

\def \id { {\mathrm{id}} }

\def \g  {\mathfrak{g}}   

\def \l {\mathfrak{l}}

\DeclareMathOperator \Span { {\mathrm{Span}} }

\newcommand \Spec { {\mathrm{Spec}} }
\begin{document}
\title[A classification of $H$-primes of quantum partial flag varieties]
{A classification of $H$-primes of quantum partial flag varieties}
\author[Milen Yakimov]{Milen Yakimov}
\address{
Department of Mathematics \\
Louisiana State Univerity \\
Baton Rouge, LA 70803 and
Department of Mathematics \\
University of California \\
Santa Barbara, CA 93106 \\
U.S.A.
}
\email{yakimov@math.lsu.edu}
\date{}
\begin{abstract}
We classify the invariant prime ideals of a quantum partial 
flag variety under the action of the related maximal torus.
As a result we construct a bijection between them 
and the torus orbits of symplectic leaves of the standard 
Poisson structure on the corresponding flag variety. 
It was previously shown by K. Goodearl and the author 
that the latter are precisely the Lusztig strata of the partial flag variety.
\end{abstract}
\maketitle
\sectionnew{Introduction}
\lb{intro}
Let $G$ be a split, simply connected, semisimple algebraic group over a field 
$\KK$ of characteristic 0 and $\g$ be its Lie algebra. Denote 
by $B$ and $B_-$ a pair of dual Borel subgroups and set $T= B \cap B_-$.
Given a set 
of simple roots $I$, one defines the standard parabolic subgroup $P_I \supset B$ 
and the multicone 
\[
\Spec \, \Big( \bigoplus_{n_i \in \Zset_{\geq 0}} 
H^0(G/P_I, \otimes_{i \notin I} \LL_{\om_i}^{n_i}) \Big)
\]
over the flag variety $G/P_I$. Here the tensor product 
involves the canonical line bundles $\LL_{\om_i}$ over 
$G/P_I$ corresponding to the fundamental weights 
for the simple roots not in $I$. Its coordinate ring 
has a canonical deformation $R_q[G/P_I]$ defined 
by Lakshmibai--Reshetikhin \cite{LR} and Soibelman \cite{S}.
The group $H$ of grouplike elements of $\UU_q(\g)$ acts 
naturally on $R_q[G/P_I]$.

Currently, little is known about the spectrum of $R_q[G/P_I]$
beyond the case of the full flag variety $G/B$. For the 
quantized ring $R_q[G/B]$ Gorelik \cite{G} classified all 
$H$-invariant prime ideals, described the inclusions between 
them and the strata of a related partition of $\Spec R_q[G/B]$.
In the general case one can apply results of Goodearl and Letzter \cite{GL}
to obtain a partition of $\Spec R_q[G/P_I]$ indexed by 
the $H$-prime ideals of $R_q[G/P_I]$, such that each stratum 
is homeomorphic to the spectrum of a Laurent polynomial ring.
The $H$-primes of $R_q[G/P_I]$ are
unknown except the case of Grassmannians which is 
due to Launois, Lenagan and Rigal \cite{LLR}.

In this note we prove a classification of the $H$-invariant 
prime ideals of the rings $R_q[G/P_I]$ associated 
to all partial flag varieties (see \thref{YY}):

\bth{main1} For an arbitrary partial flag variety $G/P_I$ the 
$H$-invariant prime ideals of $R_q[G/P_I]$  
(not containing the augmentation ideal \eqref{J+}) 
are parametrized by 
\begin{equation}
\label{SwI}
S_{W, I} : = \{ (w, v)  \in W^I \times W \mid v \leq w \}.
\end{equation}
All such ideals are completely prime.
\eth   
Here $W^I$ denotes the set of minimal length 
representatives for the elements of $W/W_I$, where $W$ is the Weyl group
of $G$ and $W_I$ is the parabolic subgroup 
of $W$ corresponding to $P_I$.

To put our results in a more geometric context, let us assume that 
the ground field is $\KK = \Cset$. The action of the torus $T$ on 
$G/P_I$ preserves the {\em{standard Poisson structure}} $\pi_I$
on $G/P_I$, cf. \cite{GY}. According to \cite[Theorem 0.4]{GY} 
its $T$-orbits of symplectic leaves are  
\[
\TT^I_{w,v} = q_I ( B \cdot w B \cap B_- \cdot v B), \; (w,v) \in S_{W,I},
\]
where $q_I \colon G/B \to G/P_I$ is the canonical projection,
cf. \cite{BGY} for the case of Grassmannians. The 
varieties $\TT^I_{w,v}$ are precisely the strata of the Lusztig 
stratification \cite{L2} of $G/P_I$ defined for the purposes of the study of 
total positivity on $G/P_I$.   

The algebra $R_q[G/P_I]$ is a quantization of the projective Poisson 
variety $(G/P_I, \pi_I)$. In particular our results provide the first step 
of the orbit method program for this situation: we obtain a bijection 
between the $H$-primes of $R_q[G/P_I]$ and the $T$-orbits of leaves of 
$(G/P_I, \pi_I)$.

The Zariski closures of $\TT^I_{w,v}$ were explicitly determined in 
\cite{GY,R}:
\[
\ol{\TT^I_{w,v}} = \{ \TT^I_{w',v'} \mid \exists z \in W_I \; \;
\mbox{such that} \; \; w \geq w' z, v \leq v' z \}.
\]
Denote by $\II_{w,v}^I$ the $H$-invariant prime ideal 
of $R_q[G/P_I]$ corresponding to $(w,v) \in S_{W,I}$ 
according to the parametrization of \thref{main1},
cf. \S \ref{3.4} for details.

Following the orbit method we make the following conjecture.

\bcj{main2} Let $(w,v), (w',v') \in S_{W,I}$, cf. \eqref{SwI}. 
Then $\II_{w,v}^I \subseteq \II_{w',v'}^I$ if and only if there exits
$z \in W_I$ such that
\[
w \geq w' z \; \; \mbox{and} \; \;  v \leq v' z.
\]
\ecj

This was established by Gorelik in the case $I= \emptyset$ of the full 
flag variety \cite{G} in which case $z$ always has to be the identity. 
In general, the conjecture is open even for Grassmannians.

Finally we prove results, which are analogous to \thref{main1}, for the 
quantum deformations \cite{LR,S} of the coordinate rings of the cones
\[
\Spec \, \Big( \bigoplus_{n \in \Zset_{\geq 0} } 
H^0( G/P_I, \LL_{ n\la} ) \Big)
\]
over $G/P_I$ for certain dominant weights $\la$.

After the completion of this paper we learned that St\'ephane Launois and
Laurent Rigal worked independently on related problems.
\\ \hfill \\
{\bf Acknowledgements.} The author is indepted to the referee for 
the careful reading of the manuscript and for numerous
detailed comments and to Ken Goodearl for very helpful correspondence.
The author was partially supported by NSF grant DMS-0701107. 
\sectionnew{Generalities on quantum groups and quantum 
flag varieties}
\lb{qalg}
\subsection{}
\label{2.1}
Let $\KK$ be a field of characteristic $0$ and $q \in \KK$ be
transcendental over $\Qset$. Let $\g$ be a split semisimple 
Lie algebra over $\KK$ of rank $r$ with Cartan matrix
$(c_{ij})$. Denote by $\UU_q(\g)$ the quantized universal 
enveloping algebra of $\g$. It is a Hopf algebra over $\KK$
with generators
\[
X^\pm_i, K_i^{\pm 1}, \; i=1, \ldots, r,
\]
as in \cite[\S 9.1]{CP}. Let $\{d_i\}_{i=1}^r$ be the standard 
choice of integers for which the matrix $(d_i c_{ij})$ 
is symmetric. Set $q_i = q^{d_i}$. Fix a nondegenerate 
invariant bilinear form $\lcor.,. \rcor$
on $\g$ such that the square length 
of a long root is equal to 2.

Let $Q$ and $Q^+$ be the sets of all integral and dominant integral 
weights of $\g$. The sets of simple roots, simple coroots, and 
fundamental weights of $\g$ will be denoted by $\{\al_i\}_{i=1}^r,$ 
$\{\al_i\spcheck\}_{i=1}^r,$ and $\{\om_i\}_{i=1}^r,$ respectively.
For $\la, \mu \in Q$ one sets $\la \geq \mu$ if 
$\mu = \la - \sum_{i=1}^r n_i \al_i\spcheck$ for some 
$n_i \in \Zset_{\geq 0}$, and $\la > \mu$ if $\la \geq \mu$ and
$\la \neq \mu$.

Recall that the weight spaces of a $\UU_q(\g)$-module $V$ are defined by
\[
V_\la = \{ v \in V \mid K_i v = q^{ \lcor \la, \al_i\spcheck \rcor} v, \; 
\forall i = 1, \ldots, r \}, \; \la \in Q.
\]
A $\UU_q(\g)$-module is a weight module if it is the 
sum of its weight spaces. The irreducible finite dimensional 
weight $\UU_q(\g)$-modules are parametrized by $Q^+$,
cf. \cite[\S 10.1]{CP} for details. For 
$\la \in Q^+$ let $V(\la)$ be the corresponding irreducible module
and $v_\la$ be a highest weight vector. All duals of finite 
dimensional $\UU_q(\g)$-modules
will be considered as left modules using the antipode of $\UU_q(\g)$.

Denote the Weyl and braid groups of $\g$ by $W$ and $\B_\g$, 
respectively. Let $s_1, \ldots, s_r$ be the simple reflections 
of $W$ and $T_1, \ldots, T_r$ be the standard generators of $\B_\g$. 
There is a natural action of $\B_\g$ on $\UU_q(\g)$ 
and the modules $V(\la)$, see \cite[\S 5.2 and \S 37.1]{L} for details.
One has $T_w ( x . v ) = (T_w x) . (T_w v)$ and 
$T_w(V(\la)_\mu) = V(\la)_{w \mu}$
for all $w \in W$, $x \in \UU_q(\g)$, $\la \in Q^+$, $v \in V(\la)$, 
$\la \in Q$.
\subsection{}
\label{2.2}
Let $G$ be the split, connected, simply connected algebraic group over $\KK$ 
with Lie algebra $\g$, and $B$ and $B_-$ be a pair of opposite Borel 
subgroups. Let $T = B \cap B_-$.

The quantized coordinate ring $R_q[G]$ is the Hopf  
subalgebra of the restricted dual of $\UU_q(\g)$ spanned by 
all matrix entries $c_{\xi, v}^\la,$ 
$\la \in Q^+$, $v \in V(\la), \xi \in V(\la)^*$:
$c_{\xi,v}^\la(x) = \lcor \xi, x v \rcor$ for $x \in \UU_q(\g)$.
One has the left and right actions of $\UU_q(\g)$ 
on $R_q[G]$:
\begin{equation}
\label{action}
x \rha c = \sum c_{(2)}(x)c_{(1)}, \; 
c \lha x = \sum c_{(1)}(x)c_{(2)}, \; 
x \in \UU_q(\g), c \in R_q[G]
\end{equation}
where $\Delta(c) = \sum c_{(1)} \otimes c_{(2)}$.

Denote by $\UU_\pm$ the subalgebras of $\UU_q(\g)$
generated by $\{X^\pm_i\}_{i=1}^r$. Let $H$ be the group 
generated by $\{K_i^{\pm 1}\}_{i=1}^r$.
The subalgebra of $R_q[G]$ invariant under the left action of $\UU_+$
will be denoted by $R^+$. It is spanned by 
all matrix entries $c_{\xi , v_\la}^\la$ where $\la \in Q^+$, 
$\xi \in V(\la)^*$ and $v_\la$ is the fixed highest weight vector 
of $V(\la)$. 

For $I \subset \{1, \ldots, r\}$ denote by $P_I \supset B$ 
the corresponding standard parabolic subgroup. Let 
$I^c = \{1, \ldots, r \} \backslash I$. Let 
$Q_I = \{ \sum_i n_i \om_i \mid i \in I^c, n_i \in \Zset\}$, 
$Q^+_I = \{ \sum_i n_i \om_i \mid i \in I^c, n_i \in \Zset_{\geq 0}\}$,
and 
$Q^{++}_I = \{ \sum_i n_i \om_i \mid i \in I^c, n_i \in \Zset_{> 0}\}$.

Denote by $\UU_q(\l_I)$ the Hopf subalgebra of $\UU_q(\g)$ 
generated by $\{X^\pm_i, K_i \}_{i \in I}$.
The quantized (multihomogeneous) coordinate ring $R_q[G/P_I]$
of the partial flag variety $G/P_I$ is defined \cite{LR,S} by
\[
R_q[G/P_I] = \Span \{ c^\la_{\xi, v_\la } 
\mid \la \in Q^+_I, \xi \in V(\la) \}.
\]
It is the subalgebra of $R_+$ invariant under the left
action \eqref{action} of the Hopf algebra $\UU_q(\l_I)$. Recall 
that each $\la \in Q_I^+$ gives rise to a 
line bundle $\LL_\la$ on the flag variety $G/P_I$.  
The ring $R_q[G/P_I]$ is a deformation of the 
coordinate ring of the multicone
\[
\Spec \, \Big( \bigoplus_{\la \in Q^+_I} H^0( G/P_I, \LL_\la) \Big)
\]
over $G/P_I$.

A subset of $R_q[G]$ is invariant under the right action of $H$ 
if and only if it is invariant under the rational action of 
the torus $(\KK^*)^r$ given by 
\[
(a_1, \ldots, a_r) \cdot c_{\xi, v_\la}^\la = 
\left( \prod_{i=1}^r a_i^{ \lcor \mu, \al_i\spcheck \rcor } \right) 
c_{\xi, v_\la}^\la, \quad
\mbox{for} \; \; \xi \in V(\la)_\mu.
\]
In particular $H$-primes and $(\KK^*)^r$-primes 
of $R_q[G/P_I]$ coincide. 
\subsection{}
\label{2.3}
Given a reduced expression
\begin{equation}
\label{wdecomp}
w = s_{i_1} \ldots s_{i_k}
\end{equation}
of an element $w \in W$, define the roots
\begin{equation}
\label{beta}
\beta_1 = \al_{i_1}, \beta_2 = s_{i_1} \al_{i_2}, 
\ldots, \beta_k = s_{i_1} \ldots s_{i_{k-1}} \al_{i_k}
\end{equation}
and the root vectors
\begin{equation}
X^{\pm}_{\beta_1} = X^{\pm}_{i_1}, 
X^{\pm}_{\beta_2} = T_{s_{i_1}} X^\pm_{i_2}, 
\ldots, X^\pm_{\beta_k} = T_{s_{i_1} \ldots s_{i_{k-1}}} X^{\pm}_{i_k},
\label{rootv}
\end{equation}
see \cite[\S 39.3]{L}. De Concini, Kac and Procesi defined 
\cite{DKP} the subalgebras $\UU_\pm^w$ of $\UU_\pm$ generated by 
$X^{\pm}_{\beta_j}$, $j=1, \ldots, k$ and proved:

\bth{DKP} (De Concini, Kac, Procesi) \cite[Proposition 2.2]{DKP}
The algebras $\UU_\pm^w$ do not depend on the choice of a 
reduced decomposition of $w$ and have the PBW basis
\begin{equation}
\label{vect}
(X^\pm_{\beta_k})^{n_k} \ldots (X^\pm_{\beta_1})^{n_1}, \; \; 
n_1, \ldots, n_k \in \Zset_{\geq 0}.
\end{equation}
\eth

The fact that the vector space spanned by the monomials 
\eqref{vect} does not depend on the choice of a
reduced decomposition of $w$ was independently obtained by Lusztig 
\cite[Proposition 40.2.1]{L}.

Recall that the universal $R$-matrix associated to $w \in W$ is 
defined by
\begin{equation}
\RR^w = \prod_{j= k, \ldots, 1} \exp_{q_{i_j}}
\left( (1-q_{i_j})^{-2} 
X^+_{\beta_j} \otimes X^-_{\beta_j} \right)
\label{Rw}
\end{equation}
where
\[
\exp_{q_i} = \sum_{n=0}^\infty q_i^{n(n+1)/2} 
\frac{n^k}{[n]_{q_i}!} \cdot
\]
In \eqref{Rw} the terms are multiplied in the order $j = k, \ldots, 1$.
The $R$-matrix $\RR^w$ belongs to a certain completion \cite[\S 4.1.1]{L}
of $\UU^w_+ \otimes \UU^w_-$ and does not depend on the 
choice of reduced decomposition of $w$. 

For all $\la \in Q^+$ and $w \in W$ fix 
$\xi_{w, \la} \in (V(\la)^*)_{- w\la}$ such that
$\lcor \xi_{w, \la}, T_w v_\la \rcor =1$. Let
\[
c^\la_w = c^\la_{\xi_{w,\la}, v_\la}.
\] 
Then $c^\la_w c^\mu_w = c^{\la+ \mu}_w= c^\mu_w c^\la_w$,
$\forall \la, \mu \in Q^+$, see e.g. \cite[\S 2.5]{Y}.

Denote
\[
c_w^I = \{ c_w^\la \mid \la \in Q^I_+ \}. 
\]

According to \cite[Lemma 9.1.10]{J} the set $c_w^{ \{1, \ldots, r \} }$ 
is Ore in $R^+$. Similarly one proves that $c_w^I$ is an Ore subset 
of $R_q[G/P_I]$.
\sectionnew{$H$-invariant prime ideals of $R_q[G/P_I]$} 
\label{hsp}
\subsection{} 
\label{3.0}
Denote by $H-\Spec_+ R_q[G/P_I]$ the set of $H$-invariant 
prime ideals of $R_q[G/P_I]$ under the right action of $H$ 
which do not contain the ideal
\begin{equation}
\label{J+}
\JJ_+^I = \Span \{ c^\la_{\xi, v_\la} \mid \la \in Q_I^{++}, \xi \in V(\la)^* \}
\end{equation} 
of $R_q[G/P_I]$.

To classify the ideals $\II$ in $H-\Spec_+ R_q[G/P_I]$,
we first partition the set according to the maximal 
quantum Schubert ideal contained in $\II$, using techniques
of Joseph \cite{J}, similar to Hodges--Levasseur \cite{HL} 
and Gorelik \cite{G}. We then relate
those strata to $H$-invariant prime ideals of the algebras $\UU^w_-$
along the lines of our previous work \cite{Y},
similarly to De Concini--Procesi \cite{DP}, and finally use 
results of M\'eriaux--Cauchon \cite{MC} and the author \cite{Y}.
\subsection{}
\label{3.1}
Recall from the introduction that $W_I$ denotes the parabolic subgroup of 
the Weyl group $W$ generated by $s_i$, $i \in I$ and $W^I$ denotes the 
set of minimal length representatives of the cosets in $W/W_I$.

We will need the following known lemma.
We include its proof for completeness.

\ble{1} Assume that $\la_j \in Q_I^+$ and $\mu_j$ are weights
of $V(\la_j)$ for $j=1,2$. Then 
$\lcor \mu_1, \mu_2 \rcor \leq \lcor \la_1, \la_2 \rcor$. 
If $\la_2 \in Q_I^{++}$, then equality implies $\mu_1 = w \la_1$
for some $w \in W^I$. If in addition $\la_1 \in Q_I^{++}$, then 
$\mu_2 = w \la_2$ for the same $w$.
\ele

\begin{proof} There exists $w \in W$ such that $w^{-1} \mu_1 \in Q^+$. 
Then $w^{-1} \mu_2 = \la_2 - \sum_{i=1}^r n_i \al_i\spcheck$, 
for some $n_i \in \Zset_{\geq 0}$ and 
\begin{multline}
\label{ineq}
\lcor \mu_1, \mu_2 \rcor = \lcor w^{-1} \mu_1, w^{-1} \mu_2 \rcor = 
\lcor w^{-1} \mu_1 , \la_2 \rcor \\ - \sum_{i=1}^r n_i \lcor w^{-1} \mu_1, \al_i \spcheck \rcor
\leq \lcor w^{-1} \mu_1 , \la_2 \rcor \leq
\lcor \la_1 , \la_2 \rcor.
\end{multline}
Assume that $\la_2 \in Q_I^{++}$ and equality holds.
Then $\la_1 - w^{-1} \mu_1 = \sum_{i=1}^r m_i \al_i\spcheck$
and $\sum_{i=1}^r m_i \lcor \la_2, \al_i\spcheck \rcor=0$. Thus $m_i = 0$
for all $i \in I^c$. Since $X_i^- v_{\la_1} = 0$ for $i \in I$ 
and $w^{-1} \mu_1 = \la_1 - \sum_{i \in I} m_i \al_i\spcheck$ is 
a weight of $V(\la_1)$, $m_i = 0$ for all $i$ and $w^{-1} \mu_1 = \la_1$.

Now assume that in addition $\la_1 \in Q_I^{++}$ and 
equality holds in \eqref{ineq}. Then $n_i = 0$ for all 
$i \in I^c$ and $w^{-1} \mu_2 = \la_2 -\sum_{i \in I} n_i \al_i\spcheck$.
Since the latter is a weight of $V(\la_2)$ in the same way we obtain 
$w^{-1} \mu_2 = \la_2$.
\end{proof}   

For $\la \in Q_I^+$ and $\mu \in Q$ denote by $\JJ_\la(\mu)$ 
the ideal of $R_q[G/P_I]$ generated by $c_{\xi, v_\la}^\la$ 
for $\xi \in V(\la)^*_{-\mu'}$, $\mu' < \mu$. 
Analogously to \cite[Proposition 9.1.5 (i)]{J} we have: 
\begin{equation}
\label{crel}
c_{\xi_1, v_{\la_1} }^{\la_1} c_{\xi_2, v_{\la_2} }^{\la_2}
- q^{ \lcor \la_1 , \la_2 \rcor - \lcor \mu_1 , \mu_2 \rcor }
c_{\xi_2, v_{\la_2} }^{\la_2} c_{\xi_1, v_{\la_1} }^{\la_1}
\in \JJ_{\la_2}(\mu_2)
\end{equation}
if $\xi_j \in V(\la_j)^*_{-\mu_j}$ for $j=1,2$. 
\subsection{} 
\label{3.2}
Following Joseph \cite[\S 9.3.8]{J} for an ideal $\II$ of $R_q[G/P_I]$
and $\la \in Q_I^+$ define 
\[
C_\II^+(\la) = \{ \mu \in Q \mid \exists \, \xi \in V(\la)^*_{-\mu}
\; \; \mbox{such that} \; \; c_{\xi, v_\la}^\la \notin \II \}.
\]
If $C_\II^+(\la)$ is empty, let $D^+_\II(\la) = \emptyset$. 
Otherwise denote by $D^+_\II(\la)$ the set of minimal elements 
of $C_\II^+(\la)$. 

\bth{2} For each prime ideal $\II$ of $R_q[G/P_I]$ which does not 
contain $\JJ_+^I$ there exists $w \in W^I$ 
such that $D_\II^+(\la) = \{ w \la\}$ for all $\la \in Q_I^+$.
\eth
For $w \in W^I$ denote by $X_w^I$ the set of those $H$-invariant 
prime ideals $\II$ of $R_q[G/P_I]$
which satisfy $D_\II^+(\la) = \{ w \la\}$ for all $\la \in Q_I^+$. 
Note that $X_w^I \subset H-\Spec_+ R_q[G/P_I]$ since $D_\II^+(\la) = \{ w \la\}$
implies that $\II$ does not contain $\JJ_+^I$. Thus we have the 
set theoretic decomposition:
\begin{equation}
\label{decomp}
H-\Spec_+(R_q[G/P_I]) = \bigsqcup_{w \in W^I} X_w^I.
\end{equation}
\\
\noindent
{\em{Proof of \thref{2}}}. We follow the idea of the proof 
of \cite[Proposition 9.3.8]{J}. Assume that $\II$ is a prime ideal 
of $R_q[G/P_I]$ which does not contain $\JJ_+^I$.

Assume that $D_\II^+(\la_j) \neq \emptyset$ and 
$\mu_j \in D_\II^+(\la_j)$ for $j = 1,2$. It follows 
from the definition of the ideal $\JJ_{\la_j}(\mu_j)$ 
(see \S 3.2) that
$\II \supset \JJ_{\la_j}(\mu_j)$ for $j=1,2$.
Fix $\xi_j \in V(\la_j)_{-\mu_j}^*$ such that 
$c^{\la_j}_{\xi, v_{\la_j}} \not\in \II$. Then 
the images $\bar{c}_j$ of $c^{\la_j}_{\xi, v_{\la_j}}$
in $R_q[G/P_I]/\II$ are normal by 
\eqref{crel} and thus are not zero divisors 
since $\II$ is prime. Applying one more time \eqref{crel} leads 
to 
\[
\bar{c}_1 \bar{c}_2 = q^{\lcor \la_1, \la_2 \rcor - 
\lcor \mu_1, \mu_2 \rcor} \bar{c}_2 \bar{c}_1 \; \;
\mbox{and} \; \;
\bar{c}_2 \bar{c}_1 = q^{\lcor \la_1, \la_2 \rcor - 
\lcor \mu_1, \mu_2 \rcor} \bar{c}_1 \bar{c}_2.
\] 
Therefore $\lcor \la_1, \la_2 \rcor= \lcor \mu_1, \mu_2 \rcor$.

Since $\II \not\supset \JJ^I_+$ there exists $\la \in Q^{++}_I$ such that 
$C_\II^+(\la) \neq \emptyset$. Then the above argument  and \leref{1} imply 
that $D_\II^+(\la) = \{ w \la \}$ for some $w \in W^I$. Let
$\xi \in V(\la)^*_{-w\la}$, $\xi \neq 0$. Since
$\dim V(\la)^*_{-w\la} = 1$, $c_{\xi, v_\la}^\la \notin \II$.
Moreover the image of $c_{\xi, v_\la}^\la$ in 
$R_q[G/P_I]/\II$ is normal
because of \eqref{crel} and all of its powers do not vanish.
As a consequence $C_\II^+(n \la) \neq \emptyset$ 
for all $n \in \Zset_{>0}$, and the above argument
and \leref{1} imply 
$D_\II^+(n \la) = \{ n w \la \}$ for the same $w$.
Furthermore, if $C^+_\II(\la')$ is nonempty 
for another $\la' \in Q^{++}_I$, then the above argument 
combined with the last assertion of \leref{1} yields  
$D^+_\II(\la') = \{ w \la'\}$ for the same $w \in W^I$.

We claim that for all $\mu \in Q^+_I$, $w \mu \in C^+_\II(\mu)$. 
There exists $n \in \Zset_{>0}$ such that $\mu \leq n \la$.
If $c_{\xi, v_\mu}^\mu \in \II$ for some $\xi \in V(\mu)^*_{-w\mu}$,
$\xi \neq 0$, then this would force 
$c_w^{n\la- \mu} c_{\xi, v_\mu}^\mu \in \II$ for 
$n \in \Zset_{>0}$ with the above property. Recall 
the well-known fact that $\dim V(\eta)^*_{-w \eta} =1$
for all $\eta \in Q_I^+$. It follows that $c^\mu_{\xi, v_\mu}$
is a nonzero scalar multiple of $c^\mu_w$. Hence, 
by \S 2.3, $c_w^{n\la- \mu} c_{\xi, v_\mu}^\mu$ 
is a nonzero scalar multiple of $c_w^{n \la}$. 
In particular, $c_w^{ n \la} \in \II$ and 
$ w (n \la) \notin C^+_\II(n \la)$, 
since $\dim V(n \la)^*_{-w(n\la)} =1$.
This is a contradiction. 

Next, we verify that $w \mu \in D^+_\II(\mu)$
for all $\mu \in Q^+_I$. If not, then
there exists $\ga = \sum_{i=1}^r n_i \al_i\spcheck$, 
$n _i \in \Zset_{\geq 0}$, $\sum_{i =1}^r n_i >0$ such that 
$w \mu - \ga \in D^+_\II(\mu)$. It follows that 
$\JJ_\la(w \mu - \ga) \subset \II$. Let 
$\xi \in V(\mu)^*_{- w \mu + \ga}$ be such that 
$c^\mu_{\xi, v_\mu} \notin \II$.  Then the image of 
$c^\mu_{\xi, v_\mu}$ in $R_q[G/P_I]/\II$ is normal
by \eqref{crel} 
and thus $c_w^\la c^\mu_{\xi, v_\mu} \notin \II$. 
This implies $w(\la+ \mu) - \ga \in C^+_\II(\la + \mu)$.
At the same time $\la + \mu \in Q^{++}_I$ forces
$D^+_\II(\la+ \mu) = \{ w(\la+\mu) \}$ which is a contradiction.

To obtain $D^+_\II(\mu) = \{ w \mu \}$ for all $\mu \in Q^+_I$,
one needs to show that 
$D^+_\II(\mu) \supsetneq \{ w \mu \}$ is impossible. If 
$\eta \in D^+_\II(\mu)\backslash \{w \mu \}$, then the argument at the beginning of
the proof and \leref{1} imply $\eta = y \mu$ for some $y \in W^I$, 
$y \neq w$. Using normality again one gets 
$w \la + y \mu \in C^+_\II(\la+ \mu)$, so 
$ w \la + y \mu \geq w(\la + \mu)$ since $\la + \mu \in Q^{++}_I$.
Thus $w \mu \leq y \mu$ which is a contradiction to $y \mu \in D^+_\II(\mu)$.
This completes the proof of the lemma.
\qed

For $w \in W^I$ define the quantum Schubert ideals
\begin{equation}
\label{Qw}
\QQ(w)_I^+ = \Span \{ c^\la_{\xi, v_\la} \mid 
\la \in Q^+_I, 
\xi \in V(\la)^*, \, \xi \perp \UU_+ T_w v_\la \}
\end{equation}
of $R_q[G/P_I]$, where ``$\perp$'' means orthogonal 
with respect to the pairing between $V(\la)^*$ 
and $V(\la)$, cf. \cite{LR,S,J,G}. The ideal $\QQ(w)_I^+$
is completely prime since it is the 
intersection of the completely prime ideal $\QQ(w)^+_{ \{1, \ldots, r\} }$
of $R^+$ (see \cite[Proposition 10.1.8]{J}) with $R_q[G/P_I]$.

\bpr{3} For all $w \in W^I$, if $\II$ is a prime ideal
of $R_q[G/P_I]$ with $D_\II^+(\la) = \{ w \la \}$ for all 
$\la \in Q^+_I$, then  
\[
\QQ(w)_I^+ \subseteq \II.
\]
\epr
\begin{proof} We use the idea of the proof of \cite[Corollary 10.1.13]{J}. 
Let $\II \in X_w^I$. We need to prove that:

(*) {\em{If $\la \in Q^+_I$, $\xi \in V(\la)_{-\mu}$, 
$\xi \perp \UU_+ v_\la$, $\mu \in Q$, then $c^\la_{\xi, v_\la} \in \II$}}.

First we show (*) for $\la \in Q^{++}_I$. Fix $\la \in Q^{++}_I$.
Assume that (*) is not correct, 
and choose $c_{\xi, v_\la}^\la$ with the property $c_{\xi, v_\la}^\la \notin \II$,
$c_{\xi, v_\la}^\la \in \QQ(w)_I^+$ such that $\mu \in Q$ is minimal.
Using $\JJ_\la(w \la) \subseteq \II$ and applying \eqref{crel}, we obtain 
\begin{equation}
\label{eeq1}
c_{\xi, v_\la}^\la c^\la_w -
q^{ \lcor \la, \la \rcor - \lcor \mu, w \la \rcor }
c^\la_w c_{\xi, v_\la}^\la \in \II.
\end{equation}
From \cite[Lemma 10.1.11 (i)]{J} one has:
\begin{equation}
\label{eeq2}
c^\la_{\xi_{1,\la}, v_\la} c_{\xi, v_\la}^\la =
q^{ \lcor \la, \la \rcor - \lcor \mu, \la \rcor }
c_{\xi, v_\la}^\la c^\la_{\xi_{1,\la}, v_\la} 
\end{equation}
(see \S 2.4 for the definition of $\xi_{y,\la}$ for $y \in W$).
Denote by $\UU^+_+$ the subalgebra of $\UU_+$ generated by 
$\{X_i^+\}_{i=1}^r$. The minimality property of $\mu$ and the fact
that $\QQ(w)^+_I$ is invariant under the
right action of $\UU_+$ (see \eqref{action}) imply that
$c^\la_{a \xi, v_\la} \in \II$ for all $a \in \UU^+_+$. Chose
$a \in \UU^+_+$ such that $\xi_{w,\la} = a \xi_{1,\la}$. 
Acting by $S(a)$ on \eqref{eeq2}, using the right action \eqref{action} 
of $\UU_+$, leads to:
\begin{equation}
\label{eeq3}
c^\la_w c_{\xi, v_\la}^\la -
q^{ \lcor \la, \la \rcor - \lcor \mu, w \la \rcor }
c_{\xi, v_\la}^\la c^\la_w \in \II.
\end{equation}
Comparing \eqref{eeq1} and \eqref{eeq3}, and using the fact 
that the image of $c^\la_w \notin \II$ in $R_q[G/P_I]/\II$
is normal implies that 
$\lcor \la, \la \rcor = \lcor \mu, w \la \rcor $.
\leref{1} implies $\mu = w \la$, which is a contraction to the 
fact that $c^\la_w \notin \QQ(w)^+_I$.

Finally we prove (*) for $\la \in Q_I^+$. Let 
$\la' \in Q_I^{++}$. If 
$c_{\xi, v_\la}^\la \in \QQ(w)^+_I$, then
$c^{\la'}_w c_{\xi, v_\la}^\la = c_{\xi', v_{\la'+\la}}^{\la'+\la}$
for some $\xi' \in V(\la'+\la)^*$, $\xi' \perp \UU_+ v_{\la'+\la}$.
Since $\la'+\la \in Q_I^{++}$, $c_{\xi', v_{\la'+\la}}^{\la'+\la} \in \II$.
Because the image of $c^\la_w \notin \II$ in $R_q[G/P_I]/\II$
is normal, $c_{\xi, v_\la}^\la \in \II$.
\end{proof}
\subsection{} 
\label{3.3}
\ble{inv} Every $\II \in H-\Spec_+ R_q[G/P_I]$ is also invariant
under the left action of $H$.
\ele
\begin{proof} Let $\II \in X_w^I$, $w \in W^I$.
For all $\la \in Q^+_I$ the images 
$\bar{c}^\la_w$ of $c^\la_w$ in $R_q[G/P_I]/\II$
are normal and do not vanish; 
thus they are not zero divisors. 
Let $\mu \in Q$, $\la_j \in Q^+_I$ and
$\xi_j \in V(\la_j)_\mu$, $j=1, \ldots, l$ are such that
$\sum_{j=1}^l c^{\la_j}_{\xi_j, v_{\la_j}} \in \II$
and $\la_1, \ldots, \la_l$ are distinct.
Eq. \eqref{crel} implies
\[
0=\big(\sum_{j=1}^l c_{\xi_j, v_{\la_j} }^{\la_j} + \II \big) \bar{c}_w^\la
= \bar{c}_w^\la 
\big( \sum_{j=1}^l 
q^{ \lcor \la_j , \la \rcor - \lcor \mu, w \la \rcor }
c_{\xi_j, v_{\la_j} }^{\la_j} + \II \big).
\]
Since $\bar{c}_w^\la$ are not zero divisors
\[
\sum_{j=1}^l 
q^{ \lcor \la_j , \la \rcor - \lcor \mu, w \la \rcor }
c_{\xi, v_{\la_j} }^{\la_j} \in \II.
\]
Since $\la \in Q_I^+$ is arbitrary, this implies that
$c_{\xi_j, v_{\la_j} }^{\la_j} \in \II$ for all $j$. So $\II$ is 
invariant under the left action of $H$.
\end{proof}

Denote by $\wt{c}^\la_w$ the image of $c^\la_w$ in $R_q[G/P_I]/\QQ(w)_I^+$.
Set $\wt{c}^I_w = \{\wt{c}^\la_w \mid \la \in Q^+_I \}$ and
\[
R_{I,w}:= \big( R_q[G/P_I]/\QQ(w)_I^+ \big) [(\wt{c}^I_w)^{-1}]. 
\]
For $\mu \in Q_I$ denote $c^\mu_w = c^{-\la_1}_w c^{\la_2}_w \in R_{I,w}$
whenever $\mu = \la_1 - \la_2$, $\la_1, \la_2 \in Q^+_I$. 
This is independent of the choice of $\la_1, \la_2$, cf. \S \ref{2.3}. 
Then: 
\[
R_{I,w}= \Span
\{ \wt{c}^{-\la'}_w (c^\la_{\xi, v_\la} + \QQ(w)^+_I)
\mid \la, \la' \in Q^+_I, \xi \in V(\la)^* \}.  
\]
Denote by $R_{I,w}^H$ the invariant subalgebra of $R_{I,w}$ with respect to the 
induced left action of $H$. We have:
\begin{equation}
\label{inv1}
R_{I,w}^H = \{ \wt{c}^{-\la}_w (c^\la_{\xi, v_\la} + \QQ(w)^+_I)
\mid \la \in Q^+_I, \xi \in V(\la)^* \}.  
\end{equation}
There is no need to take $\Span$ in the right hand side of \eqref{inv1} 
because:

(**) {\em{For all $\la, \la' \in Q^+_I$, $\xi \in V(\la)^*$, there exists
$\xi' \in V(\la + \la')^*$ such that $c_w^{-\la} c^\la_{\xi, V_\la} =
c_w^{-\la -\la'} c^\la_{\xi', v_{\la+\la'}}$.}}   

Denote by $H-\Spec R_{I,w}$ and $H-\Spec R_{I,w}^H$ the sets 
of $H$-invariant prime ideals of $R_{I,w}$ and $R_{I,w}^H$ 
with respect to the induced right action of $H$. 

If $\II \in X_w^I$, then $c^I_w \cap \II = \emptyset$ 
and $\II \supset \QQ(w)_I^+$ by \thref{2} and \prref{3}.
Therefore \cite[2.1.16(vii)]{MR} the map 
\begin{equation}
\label{mcorr}
\II \mt (\II/\QQ(w)^+_I)R_{I,w} \subset R_{I,w}
\end{equation}
defines an order preserving bijection between $X_w^I$ and $H-\Spec R_{I,w}$. 
We have $R_{I,w}= R_{I,w}^H \{\wt{c}^\mu_w \mid \mu \in Q_I \}$.
The weight lattice $Q$ acts on $R_{I,w}^H$ by ring automorphisms
by:
\begin{align*}
\mu \cdot \wt{c}^{-\la}_w (c^\la_{\xi, v_\la} + \QQ(w)^+_I) 
&= \wt{c}^{\mu}_w
\wt{c}^{-\la}_w (c^\la_{\xi, v_\la} + \QQ(w)^+_I)
\wt{c}^{-\mu}_w
\\
&= q^{\lcor \mu' - w \la, w \mu \rcor}
\wt{c}^{-\la}_w (c^\la_{\xi, v_\la} + \QQ(w)^+_I), \; \; 
\mbox{for} \;\; \xi \in V(\la)^*_{\mu'}.
\end{align*}
It is clear that
\begin{equation}
\label{sma}
R_{I,w} \cong R_{I,w}^H * Q_I
\end{equation}
(see \eqref{crel}), where $*$ stands for skew-group ring.

Let $\JJ \in H-\Spec R_{I,w}$. \leref{inv} implies that 
each ideal in $X_w^I$ is invariant under both the left and 
right actions of $H$. From the bijection \eqref{mcorr} we obtain that 
the same is true for the ideal $\JJ$ and thus
\begin{equation}
\label{invmap}
\JJ = \JJ^H \{\wt{c}^\mu_w \mid \mu \in Q_I \} =
\{\wt{c}^\mu_w \mid \mu \in Q_I \} \JJ^H, 
\end{equation}
where $\JJ^H := \JJ \cap R_{I,w}^H$. Form \eqref{sma}
one obtains that
$R_{I,w}/\JJ \cong (R_{I,w}^H/\JJ^H)*Q_I$. Since 
$R_{I,w}/\JJ$ is prime, $R_{I,w}^H/\JJ^H$ is 
$Q_I$-prime\footnote{We recall that a ring $R$ acted upon a group 
$M$ by ring automorphisms is called $M$-prime, if there are no 
nontrivial $M$-invariant ideals $\JJ_1$ and $\JJ_2$ of $R$
such that $\JJ_1 \JJ_2=0$. An $M$-invariant ideal $\II$ of $R$ is
called $M$-prime if $R/\II$ is $M$-prime.} 
(see for instance the remark after \cite[Theorem II]{P}). This
implies that $R_{I,w}^H/\JJ^H$ is $H$-prime, 
because a subset of $R_{I,w}^H/\JJ^H$ which is closed under
the right action of $H$ induced from \eqref{action} 
is necessarily closed under the action of $Q_I$.
Therefore $\JJ^H$ is an $H$-prime ideal of 
$R_{I,w}^H$. In \thref{isom} below we prove 
that $R_{I,w}^H$ and $\UU^w_-$ 
are isomorphic $H$-algebras. 
The latter is an iterated skew polynomial ring, and
Proposition 4.2 of Goodearl and Letzter \cite{GL} applies 
to give that all $H$-primes of $\UU^w_-$ are completely prime.
In particular, $\JJ^H$ is an $H$-invariant prime 
ideal of $R_{I,w}^H$.
 
In the opposite direction, 
if $\JJ_0 \in H-\Spec R_{I,w}^H$, then 
$\JJ := \JJ_0 \{\wt{c}^\mu_w \mid \mu \in Q_I \}$ is a 
two sided ideal of $R_{I,w}$ (invariant under both actions of $H$) 
and $R_{I,w}/\JJ \cong (R_{I,w}^H/\JJ_0)*Q_I$. Since 
$R_{I,w}/\JJ_0$ is prime and $Q_I$ is torsion free, 
Theorem II of Passman \cite{P} implies that 
$R_{I,w}/\JJ$ is prime and thus $\JJ \in H-\Spec R_{I,w}$.
Therefore \eqref{invmap} defines an order preserving bijection between
$H-\Spec R_{I,w}$ and $H-\Spec R_{I,w}^H$. We obtain:

\bpr{Xwmap} The map
\[
\II \in X_w^I \mt (\II/\QQ(w)^+_I)R_{I,w} \cap R_{I,w}^H 
\in H-\Spec R_{I,w}^H
\]
defines an order preserving bijection from $X_w^I$ to $H-\Spec R_{I,w}^H$.
All ideals in $X_w^I$ are completely prime.
\epr

One shows the last statement as follows. It was already indicated
that all ideals in $H-\Spec R_{I,w}^H$ are completely prime. 
If $\JJ \in H-\Spec R_{I,w}$, 
then $\JJ^H:= \JJ \cap R_{I,w}^H \in H-\Spec R_{I,w}^H$
and $R_{I,w}/\JJ \cong (R_{I,w}^H/\JJ^H)*Q_I$, thus $\JJ$ is 
completely prime. Finally, if $\II \in X_w^I$, then 
$(\II/\QQ(w)_I^+) R_{I,w} \in H-\Spec R_{I,w}$ has to be completely 
prime. Therefore $\II/\QQ(w)_I^+$ and $\II$ are completely prime too. 
\subsection{}
\label{3.4}
Similarly to \eqref{inv1} (see (**))one has:
\[
\Big( \big( R_q[G/P_I] \big)[(c^I_w)^{-1}] \Big)^H = 
\{ 
c_w^{-\la} c^\la_{\xi, v_\la} \mid 
\la \in Q^+_I, \xi \in V(\la)^* \}
\]
where the invariant subalgebra is computed with respect to 
the left action of $H$. Define
\begin{multline*}
\QQ(w)_{I,w}^+ =\{ c_w^{-\la} c^\la_{\xi, v_\la} \mid 
\la \in Q^+_I, 
\xi \in V(\la)^*, 
\xi \perp \UU_+ T_w v_\la \} 
\\
\subset \Big( \big( R_q[G/P_I] \big)[(c^I_w)^{-1}] \Big)^H,
\end{multline*}
cf. \cite[\S 6.1.2]{G} and \eqref{Qw}. Clearly $\QQ(w)_{I,w}^+$ is an ideal 
of $\Big( \big( R_q[G/P_I] \big)[(c^I_w)^{-1}] \Big)^H$ (see (**))
and one has the algebra isomorphism
\begin{multline}
\label{isom}
\Big( \big( R_q[G/P_I] \big)[(c^I_w)^{-1}] \Big)^H / \QQ(w)_{I,w}^+ \cong
R_{I,w}^H, \\
c_w^{-\la} c^\la_{\xi, v_\la} + \QQ(w)_{I,w}^+ \mt
\wt{c}_w^{-\la} \big( c^\la_{\xi, v_\la} + \QQ(w)_I^+ \big), \; \;
\la \in Q^+_I, \xi \in V(\la)^*.
\end{multline}

Analogously to the proof of \cite[Theorem 3.7]{Y}, 
cf. also \cite[Theorem 3.2]{DP}, one shows that the $\KK$-linear map 
\begin{multline}
\label{phiw2}
\phi_w \colon \Big( \big(R_q[G/P_I]\big)[(c^I_w)^{-1}] \Big)^H 
\to \UU^w_-, \\
\phi_w(c_w^{-\la} c^\la_{\xi, v_\la}) = 
(c^\la_{\xi, T_w v_\la} \otimes \id) (\RR^w), \; 
\la \in Q^+_I, \xi \in V(\la)^* 
\end{multline}
is well defined and is an $H$-equivariant algebra homomorphism.
On the first algebra one uses the right action of $H$ induced from 
\eqref{action}. On the second algebra one uses the restriction of the action 
\begin{equation}
\label{KactUw}
K \cdot x = Kx K^{-1}, \; \; K \in H, x \in \UU_q(\g)
\end{equation}
of $H$ on $\UU_q(\g)$ to $\UU^w_-$, cf. \cite[(3.18)]{Y}.

\bth{isom} The map 
$\phi_w \colon \big( \big(R_q[G/P_I]\big)[(c^I_w)^{-1}] \big)^H \to \UU^w_-$ 
is a surjective $H$-equivariant algebra homomorphism with kernel $\QQ(w)^+_{I,w}$.
It induces an $H$-equivariant algebra isomorphism between $R_{I,w}^H$ 
and $\UU^w_-$.
\eth
\begin{proof} Recall that each element of 
$\big( \big(R_q[G/P_I]\big)[(c^I_w)^{-1}] \big)^H$ is of the form
$c_w^{-\la} c^\la_{\xi, v_\la}$ for some $\la \in Q^+_I$, $\xi \in V(\la)^*$. 
It belongs to the kernel of $\phi_w$ if and only $\lcor \xi, x T_w v_\la \rcor = 0$
for all $x \in \UU^w_+$ (i.e. for all $x \in \UU_+$). This is
equivalent to $c_w^{-\la} c^\la_{\xi, v_\la} \in \QQ(w)^+_{I,w}$.

The proof of the surjectivity of $\phi_w$ is similar to the one 
of \cite[Proposition 3.6]{Y}. Assuming that $\phi_w$ is 
not surjective would imply that there exists $X \in \UU^w_+$, $X \neq 0$ 
such that $\lcor \xi, X T_w v_\la \rcor = 0$ for all $\la \in Q^+_I$,
$\xi \in V(\la)^*$. Then $X_1 = T_w^{-1} (X) \in \UU^-$, $X_1 \neq 0$ 
would satisfy $\lcor \xi, X_1 v_\la \rcor = 0$ for all $\la \in Q^+_I$,
$\xi \in V(\la)^*$. The latter is impossible since as $\UU^-$-modules
one has 
\begin{equation}
\label{ver}
V(\la) \cong \UU^- v_\la / 
\Big( \sum_{i \notin I} \UU^- 
(X_i^-)^{\lcor \la, \al_i\spcheck \rcor + 1 }v_\la \Big),
\end{equation}
see e.g. \cite[Theorem 4.3.6 (i)]{J}.

The second assertion now follows using \eqref{isom}.
\end{proof}

Recall \cite[Theorem 3.8]{Y} proved using 
Gorelik's results \cite{G}:
\bth{Y}
For the $H$-action \eqref{KactUw}, the set $H-\Spec \UU^w_-$ 
of $H$-invariant prime ideals of $\UU^w_-$ ordered under inclusion
is isomorphic to $W^{\leq w}$ as a poset. 
\eth

Eq. \eqref{decomp}, \prref{Xwmap} and \thref{Y} imply the 
main result of this note:

\bth{YY} For any quantum partial flag variety $R_q[G/P_I]$
the $H$-invariant prime ideals of $R_q[G/P_I]$ (recall \eqref{J+}) 
not containing $\JJ_+^I$ 
are parametrized by 
\begin{equation}
\label{doubleset}
\{ (w, v)  \in W^I \times W \mid v \leq w \}.
\end{equation}
All such ideals are completely prime.
\eth
Denote by $\II^I_{w,v}$ the $H$-invariant prime ideal of $R_q[G/P_I]$
in $X_w^I$ which corresponds to the ideal $I_w(v)$ of \cite[Theorem 3.8]{Y}
under the order preserving bijection from \prref{Xwmap} and 
the isomorphism from \thref{isom}. Tracing back
those bijections and using the poset part of the statement 
of \thref{Y}, one obtains that for all $v, v' \leq w$:
\begin{equation}
\label{easyCor}
\II^I_{w, v} \subseteq \II^I_{w, v'} \; \; \mbox{if and only if}
\; \; v \leq v'.
\end{equation}
This proves the special case of \cjref{main2} when $w = w'$, 
but the general statement is harder. 

\bre{notroot} One can define the algebras $\UU_q(\g)$, $R_q[G/P_I]$, 
$\UU^w_-$ over any field $\KK$ (not necessarily of characteristic 0), 
for $q \in \KK$ which is not a root of unity. In this more general 
setting M\'eriaux and Cauchon \cite{MC} proved that the $H$-invariant 
prime ideals of $\UU^w_-$ are parametrized by $W^{\leq w}$ (though 
the inclusions between them are unknown). All results of this section
trivially carry out to this more general setting. As a result one obtains 
that there is a bijection between 
$H - \Spec_+ R_q[G/P_I]$ and the set \eqref{doubleset} 
for the case when $\UU_q(\g)$ is defined over an arbitrary field 
$\KK$ and $q \in \KK$ is not a root of unity.   
\ere
\subsection{}
\label{3.5}
Throughout this subsection fix $\la \in Q^{++}_I$. Consider the subalgebra
\cite{LR,S}:
\[
R_q[G/P_I]^\la = \Span \{ c^{n \la}_{\xi, v_{n \la} } 
\mid n \in \Zset_{\geq 0},  \xi \in V( n \la) \}
\]
of $R_q[G/P_I]$. It is a deformation of the 
coordinate ring of the cone
\[
\Spec \, \Big( \bigoplus_{n \in \Zset_{\geq 0} } H^0( G/P_I, \LL_{ n\la} ) \Big)
\]
over $G/P_I$ associated to $\la \in Q^{++}_I$, cf. \S \ref{2.2} for the 
definition of the line bundles $\LL_{n \la}$.

Define the ideal 
\[
\JJ_+^\la = \Span \{ c^{n \la}_{\xi, v_{n \la} } \mid n \in \Zset_{>0}, 
\xi \in V(n \la)^* \}
\]
of $R_q[G/P_I]^\la$, cf. \eqref{J+}.  
Denote by $H-\Spec_+ R_q[G/P_I]^\la$ the set of $H$-invariant 
prime ideals of $R_q[G/P_I]^\la$ under the right action \eqref{action}
of $H$ which do not contain the ideal $\JJ_+^\la$.

For an ideal $\II$ of $R_q[G/P_I]^\la$
and $n \in \Zset_{>0}$ define 
\[
C_\II^+(n) = \{ \mu \in Q \mid \exists \, \xi \in V(n \la)^*_{-\mu}
\; \; \mbox{such that} \; \; c_{\xi, v_{ n \la}}^{ n \la} \notin \II \}.
\]
If $C_\II^+(n) = \emptyset$, let $D^+_\II(n) = \emptyset$. 
Otherwise denote by $D^+_\II(n)$ the set of minimal elements 
of $C_\II^+(n)$. Denote the quantum Schubert ideal \cite{LR,S,J,G}: 
\[
\QQ(w)_\la^+ = \Span \{ c^{n \la}_{\xi, v_{n \la} } \mid 
n \in \Zset_{>0}, \xi \in V(\la)^*, \, \xi \perp \UU_+ T_w v_\la \}
\]
of $R_q[G/P_I]^\la$, cf. \eqref{Qw}. Analogously to \thref{2} and
\prref{3} one shows:
 
\bpr{la1} (1) For each prime ideals $\II$ of $R_q[G/P_I]^\la$ which does not 
contain $\JJ_+^\la$ there exists $w \in W^I$ 
such that $D_\II^+(n) = \{ n w \la\}$ for all $n \in \Zset_{>0}$.

(2) For a given $w \in W^I$, all prime ideals of $R_q[G/P_I]^\la$
satisfying the condition in (1) contain the ideal $\QQ(w)^+_{\la}$. 
\epr

Given $w \in W^I$, let $X_w^\la$ be the set of $H$-invariant prime ideals 
$\II$ of $R_q[G/P_I]^\la$ such that $D_\II^+(n) = \{ n w \la \}$, 
$\forall n \in \Zset_{>0}$. Then $X_w^\la \subset H-\Spec_+ R_q[G/P_I]^\la$
and
\[
H-\Spec_+ R_q[G/P_I]^\la = \bigsqcup_{w \in W^I} X_w^\la.
\]
Similarly to \cite[Lemma 9.1.10]{J} one shows that 
$\{(c_w^\la)^n\}_{n \in \Zset_{\geq 0}}$ is an Ore subset of $R_q[G/P_I]^\la$.
Let $\ol{c}^{\la}_w$ be the image of $c^{\la}_w$ in $R_q[G/P_I]^\la/\QQ(w)_\la^+$.
Set
\[
R_{\la,w}:=\big( R_q[G/P_I]^\la/\QQ(w)_\la^+ \big) [(\ol{c}^{\la}_w)^{-1}]. 
\]
Consider the induced left action of $H$ on $R_{\la, w}$ from \eqref{action}, 
and denote by $R_{\la,w}^H$ the corresponding invariant subalgebra. Similarly 
to \prref{Xwmap} one establishes that:

{\em{There is an order preserving bijection between 
$X_w^\la$ and $H-\Spec R_{\la, w}^H$ given by
\begin{equation}
\label{bij}
\II \in X_w^\la \mt (\II/\QQ(w)^+_\la)R_{\la,w} \cap R_{\la,w}^H 
\in H-\Spec R_{\la,w}^H,
\end{equation}
where $H-\Spec$ refers to the set of $H$-invariant prime ideals 
with respect to the induced right action from \eqref{action}. 
All ideals in $X_w^\la$ are completely prime.
}}

\noindent
One has (see (**)):
\[
\Big( \big( R_q [G/P_I]^\la \big)[(c^\la_w)^{-1}] \Big)^H = 
\{ 
c_w^{-n \la} c^{n \la}_{\xi, v_{n \la}} \mid 
n \in \Zset_{\geq 0}, \xi \in V(n \la)^* \}
\]
where $(.)^H$ denotes the invariant subalgebra with respect to the 
induced left $H$-action from \eqref{action}. 

Similarly to \thref{isom} one proves:

\bpr{phila} The $\KK$-linear map 
\begin{multline*}
\psi_w \colon \Big( \big(R_q[G/P_I]^\la\big)[(c^\la_w)^{-1}] \Big)^H 
\to \UU^w_-, \\
\psi_w(c_w^{-n \la} c^{n \la}_{\xi, v_{n \la} }) = 
(c^{n \la}_{\xi, T_w v_{n \la} } \otimes \id) (\RR^w), \; 
n \in \Zset_{\geq 0}, \xi \in V(n \la)^* 
\end{multline*}
is an $H$-equivariant surjective algebra homomorphism with kernel
\[
\{ c_w^{-n \la} c^{ n \la}_{\xi, v_{ n \la} } \mid n \in \Zset_{> 0}, 
\xi \in V(n \la)^*, \xi \perp \UU_+ T_w v_{n \la} \}.
\]
Here $H$ acts on the the first algebra by the induced right action from 
\eqref{action} and on the second algebra by \eqref{KactUw}.

The homomorphism $\psi_w$ induces an $H$-equivariant algebra isomorphism between 
$R_{\la,w}^H$ and $\UU^w_-$.
\epr
Invoking \thref{Y}, one obtains:

\bth{laHsp} For all $\la \in Q^{++}_I$, the
$H$-invariant prime ideals of $R_q[G/P_I]^\la$ not containg 
$\JJ_+^\la$ are parametrized by the set
\[
\{ (w, v)  \in W^I \times W \mid v \leq w \}.
\]
All such ideals are completely prime.
\eth

Denote by $\II^\la_{w, v}$ the ideal of $R_q[G/P_I]^\la$
which corresponds to the ideal $I_w(v)$ of $\UU^w_-$ of 
\cite[Theorem 3.8]{Y} under the bijections of
\eqref{bij} and \prref{phila}. We conjecture:

\bcj{main3} Let $(w,v), (w',v') \in S_{W,I}$, cf. \eqref{SwI}.
One has
$\II_{w,v}^\la \subseteq \II_{w',v'}^\la$ if and only if there exits
$z \in W_I$ such that
\[
w \geq w' z \; \; \mbox{and} \; \;  v \leq v' z.
\]
\ecj
Analogously to \eqref{easyCor} one uses the poset part of the statement of 
\thref{Y}, the order preserving bijections \eqref{bij} and
\prref{phila} to prove the case of \cjref{main3} when $w = w'$. 

\bre{notrootb} Similarly to \reref{notroot} one can define the algebras 
$R_q[G/P_I]^\la$, over any field $\KK$ (not necessarily of characteristic 0), 
for $q \in \KK$ which is not a root of unity. The above arguments
and the M\'eriaux--Cauchon \cite{MC} result parametrizing 
$H$-invariant prime ideals of $\UU^w_-$ prove that the 
parametrization of $H$-primes of $R_q[G/P_I]^\la$ from
\thref{laHsp} is valid in this more general situation.
\ere

\end{document}